\documentclass[12pt,leqno]{article}
\usepackage{amsfonts,amsmath,amssymb,epsfig,fancyhdr,graphicx,mathrsfs,eucal,yfonts,dsfont}
\usepackage{flafter}
\usepackage{mathrsfs}
\usepackage{natbib}
\usepackage{url}

\def\BibTeX{{\rm B\kern-.05em{\sc i\kern-.025em b}\kern-.08em
    T\kern-.1667em\lower.7ex\hbox{E}\kern-.125emX}}

\numberwithin{equation}{section}

\newtheorem{theorem}{Theorem}

\newtheorem{corollary}[theorem]{Corollary}

\newtheorem{definition}[theorem]{Definition}

\newtheorem{lemma}[theorem]{Lemma}

\newenvironment{proof}[1][Proof]{\textbf{#1.} }{\ \rule{0.5em}{0.5em}}

\newcommand\aff{\operatorname{\mathsf{aff}}}
\newcommand\Paths{\operatorname{Paths}}

\newcommand\conv{\operatorname{\textsf{conv}}}
\newcommand\cone{\operatorname{\textsf{cone}}}

\renewcommand\emptyset{\varnothing}

\newcommand\rinter{\operatorname{\textsf{rint}}}

\newcommand\rel{\operatorname{rel}}

\newcommand{\vL}{\overrightarrow{L}}

\renewcommand{\phi}{\varphi}
\renewcommand{\epsilon}{\varepsilon}

\newcommand{\0}{\mbox{\boldmath$0$}}
\newcommand{\bb}{\mbox{\boldmath$b$}}

\newcommand{\bl}{\mbox{\boldmath$l$}}

\newcommand{\p}{\mbox{\boldmath$p$}}

\newcommand{\q}{\mbox{\boldmath$q$}}

\newcommand{\x}{\mbox{\boldmath$x$}}

\newcommand\cSH{\mathcal{S}H}

\newcommand\ff{f}

\newcommand\bH{\mathbf{H}}

\newcommand\bP{\mathbf{P}}

\newcommand\bT{\mathbf{T}}

\newcommand\N{\mathbb{N}}
\newcommand\R{\mathbb{R}}
\newcommand\HH{\mathbb{H}}
\newcommand\SSS{\mathbb{S}}
\newcommand\X{\mathbb{X}}

\newcommand\cC{\mathcal{C}}
\newcommand\cD{\mathcal{D}}
\newcommand\cG{\mathcal{G}}
\newcommand\cF{\mathcal{F}}
\newcommand\cI{\mathcal{I}}

\newcommand\cS{\mathcal{S}}
\newcommand\cT{\mathcal{T}}

\newcommand\cL{\mathscr{L}}
\newcommand\cM{\mathscr{M}}
\newcommand\cN{\mathscr{N}}
\newcommand\cU{\mathscr{U}}

\newcommand\mfH{\textfrak{H}}

\begin{document}

\title
{Convexity of Hypersurfaces in Spherical  Spaces}
\author{Konstantin Rybnikov\footnote{Email: \protect  \url{Konstantin_Rybnikov@uml.edu}}}
\date{\today}
\vspace{-1in}
\maketitle
\begin{abstract}
A spherical set is called convex if for every pair of its points there is at
least one minimal geodesic segment that joins these points and lies in the set.
We prove that for $n \ge 3$ a complete locally-convex (topological) immersion of
a connected $(n-1)$-manifold into the $n$-sphere is a surjection onto the
boundary of a convex set.
\end{abstract}

\section{Introduction}

Van Heijenoort (1952) proved that a complete locally-convex immersion $\ff$ of a
connected manifold $\cM$ ($\dim~\cM=n-1$) into $\R^n$ ($n\ge3$) is a
homeomorphism onto the boundary of a convex body, provided $\ff$ has a point of
strict convexity. For $n=3$ this result, according to Van Heijenoort himself,
follows from four theorems in A.D. Alexandrov's book (1948).

Suppose  now $f:\cM \rightarrow \R^n$ does not have a point of strict convexity.
Jonker \& Norman (1973) proved that when $f: \cM \rightarrow \R^n$ is \emph{not}
a homeomorphism onto the boundary of a convex body, $\ff(\cM)$ is the direct
affine product of a locally-convex plane curve and a complementary subspace $L
\cong \R^{n-2}$ of $\R^n$. On the other hand, if $f$ is still a homeomorphism
onto the boundary of a convex body, they showed that $\ff(\cM)$ is the direct
product of a compact convex hypersurface in a $(g+1)$-subspace ($g \ge 0$) and a
complementary subspace $L \cong \R^{n-g-1}$ of $\R^n$.

The question of sufficient conditions for convexity of a hypersurface in $\SSS^n$ naturally appears in the important computational geometry problem of checking convexity of a PL-hypersurface in Euclidean space (e.g. see Rybnikov, 200X). To observe this connection just notice that the shape of a PL-hypersurface  (in $\R^n$) at a vertex $v$ is described by a hypersurface in $\SSS^{n-1}$ obtained by intersecting the star of the vertex with a small sphere centered at $v$.   Convexity checkers for surfaces in $\R^2$ and $\R^3$ have been implemented in the LEDA system for computational geometry and graph theory (Mehlhorn \& N\"{a}her, 2000).
We show that for locally-convex immersions into a sphere of dimension $n \ge 3$
the absence of  points of strict convexity  cannot result in the loss of global
convexity, as it happens in the Euclidean case:

\begin{theorem}\label{theorem:Spherical_Lemma}
Let $i: \cM \rightarrow \SSS^n$ ($n\ge3$) be a complete locally-convex immersion
of a connected $(n-1)$-manifold $\cM$. Then
$i(\cM) = \mathbb{S}^n \cap \partial K$, where $K$ is a convex cone in
$\mathbb{R}^{n+1}$ containing the origin.
\end{theorem}
The proof of this theorem
relies on the result of Jonker \& Norman, although their theorem does not
directly imply ours. One of the difficulties is that a compact convex set on the
sphere may be free of extreme points! We can observe the following ``tradeoff" between
convexity and bijectivity requirements for complete hypersurfaces without boundary immersed in $\X^n$, where $\X^n$ is one of $\R^n$, $\SSS^n$, $\HH^n$ for $n \ge 3$. In Euclidean space a locally-convex hypersurface may fail to be convex, but if it is convex, the surface is always embedded. In the spherical space a locally-convex hypersurface always bounds a convex body, but the surface need not be embedded.
It is worth noting that in the hyperbolic space a locally-convex hypersurface need not be neither convex, nor embedded.
\section{Notions and notation}\label{sec:definitions}  From now on
$\X$ (or $\X^n$) denotes $\R^n$ or $\SSS^n$, where $n \in \N=\{0,1,\dots\}$. All maps are
continuous.
\begin{definition}
A \emph{surface in } $\X$ is a pair $(\cM,r)$ where $\cM$ is a manifold,
 with or without boundary, and $r:\cM \rightarrow \X$ is a continuous map,
hereafter referred to as the \emph{realization} map.
\end{definition}

Let $\dim \cM+1=\dim \X$. Then  $(\cM,r)$ is called \emph{locally convex at} $p \in \cM$ if
we can find a neighborhood of $p$,  $\cN_p \subset \cM$, and a convex body $K_p \subset \X$
for $p$ such that $r|_{\cN_p}:\cN_p \rightarrow r(\cN_p)$ is a homeomorphism and $r(\cN_p)
\subset \partial K_{p}$. This definition was introduced by Van Heijenoort. A.D.
Alexandrov's (1948) concept of local convexity, although restricted to more specific
classes of surfaces, is essentially equivalent to Van Heijenoort's one. We refer to $K_p$
as a convex witness for $p$. (Here, as everywhere else, the subscript indicates that $K_p$
depends on $p$ in some way but is not necessarily determined by $p$ uniquely.) Thus, the
local convexity at $\p=r(p)$ may fail because $r$ is not a local homeomorphism at $p$ or
because no neighborhood $\cN_p$ is mapped by $r$ into the boundary of a convex body, or for
both of these reasons. In the first case we say that the immersion assumption is violated,
while in the second case we say that the convexity is violated. Often, when it is clear
from the context that we are discussing the properties of $r$ near $\p=r(p)$, we  say that
$r$ is convex at $\p$.  If $K_{p}$ can be chosen so that $K_{p} \setminus r(p)$ lies in an
open half-space defined by some hyperplane passing through $r(p)$, the realization $r$ is
called \emph{strictly convex} at $p$. We will also sometimes refer to $(\cM,r)$ as strictly
convex at $r(p)$. We will often apply local techniques of Euclidean convex geometry to
$\SSS^n$ without restating them explicitly for the spherical case.

The relative boundary $\rel \partial C$ of a convex set $C$ is a manifold. In the above
definition of local convexity convex bodies can be replaced with convex sets and boundaries
of convex bodies with relative boundaries of convex sets.   Such a modified definition,
which is equivalent to the traditional one in the context of our paper, has the following
advantage. Without specifying dimensions of $\cM$ and $\X$, we can say that $r:\cM
\rightarrow \X$ is locally convex at $x \in \cM$ if there is a neighborhood $\cN_x$ such
that $r|_{\cN_x}$ is a homeomorphism onto a neighborhood of $r(x)$ in $\rel
\partial K_x$, where $K_x$ is a convex set in $\X$. Such a relativized definition would make notation and
formulations more concise and aesthetically pleasing; however, to be consistent with
previous works on the subject, we follow the traditional definition.

To avoid a common confusion caused by (at least three) different usages of  \verb"closed" in
English texts on the geometry-in-the-large,  we use this word for closed subsets of topological spaces
only. We will not use the term ``closed surface" at all; a closed submanifold stands for a submanifold
 which happens to be a closed subset in the ambient manifold.  Whenever we want to include manifolds
 with boundary  into our considerations we explicitly say so.

 A map $i:\cM \rightarrow \X$ is called an \emph{ immersion} if $i$ is a
local homeomorphism; in such a case we also refer to $(\cM,i)$ as a surface
immersed into $\X$. This is a common definition of immersion in the context of
non-smooth geometry in the large (e.g. see Van Heijenoort, 1952); a more restrictive
definition is used in differential geometry and topology, furthermore, some authors
define an immersion as a continuous local bijection. Although the latter definition
 is not, in general, equivalent to the common one, it is equivalent to the common one
 in the context of the theorems stated in this paper.
 A map $e:\cM \rightarrow \X$ is called an
\emph{embedding} if $e$ is a homeomorphism onto $e(\cM)$. Obviously, an embedding is an
immersion, but not vice versa.

A set $K \subset \X$ is called \emph{convex} if for any $x,y \in K$ there is a
geodesic segment of minimal length with end-points $x$ and $y$ that lies in $K$.
Right away we conclude that the empty set and all one point sets are convex.
 A \emph{convex body} in $\X$ is a closed convex set of full dimension; \emph{a convex body may be unbounded.} Let $\dim \cM+1 = \dim \X$. Then 
 $(\cM,r)$ is called \emph{locally convex at} $p \in \cM$ if we can find
a neighborhood $\cN_p \subset \cM$ of $p$ and a convex body $K_p \subset \X$  such that $r|_{\cN_p}:\cN_p \rightarrow r(\cN_p)$ is a homeomorphism and
$r(\cN_p) \subset \partial K_{p}$. In such a case we refer to $K_p$ as a convex witness
for $p$. (Here, as everywhere else, the subscript indicates that $K_p$ depends on $p$ in some way but is not necessarily determined by $p$ uniquely.) Thus, the local convexity at $\p=r(p)$ may
fail because $r$ is not a local homeomorphism at $p$ or because no neighborhood
$\cN_p$ is mapped by $r$ into the boundary of a convex body, or for both of these
reasons. In the first case we say that the immersion assumption is violated,
while in the second case we say that the convexity is violated. Often, when it is
clear from the context that we are discussing the properties of $r$ near
$\p=r(p)$, we  say that $r$ is convex at $\p$.  If $K_{p}$ can be chosen so that
$K_{p} \setminus r(p)$ lies in an open half-space defined by some hyperplane
passing through $r(p)$, the realization $r$ is called \emph{strictly convex} at
$p$. We will also sometimes refer to $(\cM,r)$ as strictly convex at $r(p)$.
We will often apply local techniques of Euclidean
convex geometry to $\SSS^n$ without restating them explicitly for the spherical
case.

We gave a metric definition of convexity in $\SSS^n$. One can argue that convexity is
intrinsically affine notion. In differential geometry the geodesic property can be defined
locally via the notion of affine connection $\delta$. A (directed) geodesic segment is then defined
as a smooth curve $\gamma$ from $x$ to $x'$ such that

  \[\nabla_{\dot\gamma(t)}\dot\gamma(t) = 0\; \textrm{for all}\; t \in [0,1], \]

and the curve is minimal with respect to containment, i.e. $x$ and $x'$ occur only as
the source and the of this curve. Since there is no metric in a space with affine
connection, it is meaningless to talk about the shortest geodesics. Under this approach a
convex set is a set where every two points can be connected by \emph{a} geodesic segment
lying in the set. From this prospective  the segment obtained by going \emph{clockwise}
from $0$ to $\pi/2$ on $\SSS^1$ is just as good as the segment $[0,\pi/2]$. The concept of
locally-convex hypersurface is free of concepts of distance and orientation: e.g. we do not
say from which side the convex witness abuts the surface. The fact that a locally-convex
hypersurface in $\R^n$ is orientable is a \emph{consequence} of local convexity.

Certain general convexity results can be proven in the context of a manifold with affine
connection, and even in more general spaces, that include surfaces of (not necessarily
convex) polyhedra and ball-polyhedra in $\R^n$. For example, Klee's theorem that the
boundary of a convex set is a disjoint union of convex faces of various dimensions (e.g. Rockafellar, 1997, p. 164) can be proven for locally-convex hypersurfaces in very general spaces with geodesics (Rybnikov, 2005). Here by a \emph{space with geodesics} we mean a much more general space than spaces under the same name considered by Busemann and Gromov (their notions are metric -- see e.g. Berger, 2003, pp. 678-680). The notion of geodesic and convexity in such a space
should be defined locally in a sheaf-theoretic fashion. 
We will not give any formal
definitions here, but only indicate, informally, that such spaces are, although being
$C^0$ ``in the worst case" are essentially piecewise-analytic geometries which are nice
enough to avoid bizarre topological behavior.  
 The philosophy here is that any good definition must be, in an appropriate
 sense, of local nature. Hence, the requirement that a geodesic segment must be globally of minimal length seems premature in the general study of geometric convexity.
We will now show that for $\SSS^n$ the affine definition of convexity can be used in the
context of the theory of locally-convex hypersurfaces considered in this paper. A set $S$
in $\SSS^n$ is called $A$-convex if for any $x,x' \subset S$ there is a geodesic segment
(in the sense of affine convexity -- see above) that joins $x$ and $x'$ and lies in $S$. Below we consider $\SSS^n$ centered at $\0$ and embedded into $\R^{n+1}$. $\cone (S)$ stands for $\{\R_+x\;\vline\; x \in S\}$.
\begin{lemma}
If $S \subset \SSS^n$ is convex and does not belong to any subspace of dimension less than $n$, then $\overline{\SSS^n \setminus S}$ is $A$-convex.
\end{lemma}
\begin{proof} Let $x,x' \in \overline{\SSS^n \setminus S}$. Consider $K=\aff\{\0,x,x'\} \cap \cone(S)$. Since both of these sets are convex cones with common apex, their intersection $K$  is also a  convex cone with apex at $\0$.

If $K$ is a linear subspace of $\R^{n+1}$, then there is supporting hyperplane $H$ through $\aff\{\0,x,x'\}$ for $\cone(S)$.  In this case $\cone(S)$ lies in one of the closed halfspaces defined by $H$. Thus $H \cap \SSS^n \subset \overline{\SSS^n \setminus S}$,  which means that any arc of a great circle which passes through $x$ and $x'$ and lies in $H$ belongs to $\overline{\SSS^n \setminus S}$.

If $K$ is not a subspace of $\R^{n+1}$, then one of the arcs of $\aff K \cap \SSS^n$ with endpoints $x$ and $x'$ lies outside of $S$.
\end{proof}
\medskip
  Let us recall (see e.g. Rockafellar, 1997) that a point $\p$
on the boundary of a convex set $C$ is called \emph{exposed } if $C$ has a support
hyperplane that intersects $\overline{C}$,  the closure of $C$, only at $\p$.
Thus, an \emph{exposed} point on a convex body $K$ is a \emph{point of strict
convexity} on the hypersurface $\partial K$. Conversely, for a point of strict convexity
$p\in \cM$ for $(\cM,r)$ the image $i(p)$ is an exposed point of any convex witness for $p$. Local convexity can be defined in many other, non-equivalent, ways (e.g., see van Heijenoort).

A hypersurface $(\cM,r)$ is (globally) \emph{convex} if there exists a convex body $K
\subset \X^n$ such that $r$ is a homeomorphism onto $\partial K$. Hence, we exclude the
cases where $r(\cM)$ is the boundary of a convex body, but $r$ fails to be injective.
Of course, the algorithmic and
topological aspects of this case may be interesting to certain areas of  geometry,
e.g. origami.

\section{Geometry of locally-convex immersions}\label{sec:locally sphere}
Recall that a path joining points
$x$ and $y$ in a topological space ${\cal T}$ is a map $\alpha:[0,1]
\rightarrow {\cal T}$, where $\alpha(0)=x$ and $\alpha(1)=y$. Denote by $\Paths_{\cM}(x,y)$ the set of all paths
joining $x,y \in \cM$.

Any realization $r: \cM \rightarrow \X^n$ induces a distance $d_{r}$ on $\cM$ by

\[ d_r(x,y)=\underset{\alpha \in {\Paths_{\cM}(x,y)}}{\inf} { |r(\alpha)| }, \]

\par \noindent where $|r(\alpha)| \in \R_+ \cup \infty$ stands  for the length of the
$r$-image of the path $\alpha$ joining $x$ and $y$ on $\cM$ (we call it the
$r$-\emph{distance}, because it is not always a metric).

 Of course, for a general
 realization $r$ it is not clear \emph{a priori} that there is a path  of
finite length on $r(\cM)$ joining $r(x)$ and $r(y)$, even when  $x$ and $y$ are in the same connected component.
The notion of \emph{complete} realization is essential to the correctness of van Heijenoort's theorem. A realization $r:\cM \rightarrow \X$ is called \emph{complete} if every Cauchy sequence on $\cM$ (with respect to the distance induced by $r$ on $\cM$) converges. Completeness is a rather subtle notion: a space may be
complete under a metric $d$ and not complete under another metric $d_1$, which is
 topologically equivalent to $d$ (i.e. $x_n \overset{d}{\rightarrow} a$ iff $x_n \overset{d_1}{\rightarrow} a$).

A realization is called \emph{proper} if the preimage of every compact set is compact. A proper realization is always closed. For any given natural class of realizations (e.g. PL-surfaces, semialgebraic surfaces, etc) it is usually much easier to check for properness than for completeness. Furthermore, the notion of properness is topological, while that of completeness is metrical.  Note that in some sources, such
as the paper by (Burago and Shefel, 1992), completeness with respect to the $r$-metric is
called \emph{intrinsic completeness,} while properness is referred to as
\emph{extrinsic completeness.} The following is well-known for immersions (see e.g. Burago and Shefel, p. 50), but is also true for arbitrary proper realizations. The proof given here was suggested by Frank Morgan.

\begin{lemma}\label{lem:proper complete} A proper realization $r$ of any  manifold $\cM$ in $\X$ is complete.
\end{lemma}
\begin{proof} Let $\{x_n\} \subset \cM$ be Cauchy. Then $\{r(x_n)\}$ is also Cauchy in the $r$-distance and, therefore, in the intrinsic distance of $\X$ as well. Since $\X$ is complete,
$\{r(x_n)\}$ converges to some point $y$ of $\X$. Since $r(\cM)$ is closed, $y \in r(\cM)$.

For any $k \in \N$ there is $j(k)$ such that for any $i \ge j(k)$ we have $d_r(x_i,x_{i+1})<\frac{1}{2^{k}}$. Note that in this case $\sum_{k\in \N} d_r(x_{j(k)},x_{j(k+1)})$ converges. As  $\{r(x_n)\}$ is convergent, it lies in some compact set $S \subset \X$.  Since  $r$ is proper, $r^{-1}S$ is compact. Thus, $x_n$ have an accumulation point $x$. As $r$ is continuous and $r(x_n) \rightarrow y$ in $\X$, $r(x)=y$.

Let us show that $x_{j(k)}$ converges to $x$ in the $r$-distance.
For each $k$  there is a path $p_k$ of length less than $\frac{1}{2^{k}}$ (in the $r$-distance)  from $x_{j(k)}$ to $x_{j(k+1)}$. For each $k$ we can form a path $\alpha_k$ with source  $x_{j(k)}$ by concatenating $p_i$,  $p_{i+1}$,...,etc, for all $i \ge k$. Since  $\{x_{j(k)}\}$ converges to $x$, $\alpha_k$ is a path from $x_{j(k)}$ to $x$. Since  $\sum_{k\in \N} d_r(x_{j(k)},x_{j(k+1)})$ converges, it is a path of finite length. Thus, $\{x_{j(k)}\}$ converges to $x$ in the $r$-distance. Since a subsequence of $\{x_n\}$ converges to $x$ in the $r$-distance, $\{x_n\}$ also converges to $x$ in the $r$-distance.
\end{proof}

The reverse implication is true for locally-convex
immersions, but not, for example, for saddle surfaces (e.g. see
Burago and Shefel, p. 50):

\begin{lemma} (Van Heijenoort) \label{lem:complete proper} A complete locally-convex immersion of a connected $(n-1)$-manifold into $\X^n$ is proper.
\end{lemma}

\begin{lemma}\label{lemma:arcwise}(Van Heijenoort, 1952; pp. 227-228)
Let $f:\cM \rightarrow \X^n$  be a complete locally-convex immersion of an
$(n-1)$-manifold $\cM$. Then any two points in the same connected component of $\cM$ can be connected by an arc of finite length. The topology on $\cM$ defined by the $f$-distance is equivalent to the intrinsic (original) topology on $\cM$.
\end{lemma}

Van Heijenoort's proofs of  Lemmas \ref{lemma:arcwise} and \ref{lem:complete proper} given in the original for
 $\R^n$  are valid, word by word, for $\SSS^n$ and $\mathbb{H}^n$, since these lemmas are entirely of local nature.
If $f$ is a locally-convex immersion, then for a ``sufficiently small" subset $\cS$ of $\cM$ the map $f \vert_{\cS}$ is a homeomorphism and, therefore, the topology on $\cS$ that is induced by the metric topology of $\X^n$ is equivalent to the intrinsic topology of $\cS$ and, thanks to Lemma \ref{lemma:arcwise}, to the $f$-distance topology. Thus, when $f$ is a complete locally-convex immersion, then for sufficiently small subsets of $\cM$ (but not $i(\cM)$ !) the three topologies considered in this section are equivalent -- this fact will be used throughout the text  without an explicit reference to the above lemmas. The following is our starting point.

\begin{theorem} \emph{(Van Heijenoort, 1952)}\label{theorem:H}
If a complete locally-convex immersion $\ff$ of a connected $(n-1)$-manifold $\cM$ into
$\R^n$ ($n\ge3$) has a point of strict convexity, then $\ff$ is a homeomorphism onto
the boundary of a convex body.
\end{theorem}
There is no need to check the existence of a point of strictly convexity in the compact
case.

\begin{lemma}\label{lemma:local_is_enough} If $f:\cM \rightarrow \R^n$ is a locally-convex immersion of a compact
connected $(n-1)$-manifold $\cM$, then $f$ has a point of strict convexity.
\end{lemma}
\begin{proof}As $\cM$ is compact and $f$ is an immersion,
$\conv f(\cM)$ is a compact subset of $\R^n$.  Since $\conv f(\cM)$ is compact, it is
also bounded and, in particular, does not contain lines. Any non-empty convex set,
which is free of lines, has a non-empty set of extreme points (a point on the boundary
of a convex set is extreme if it is  not interior to any line segment contained in
set's boundary). Thus $\partial \conv f(\cM)$
contains an \emph{extreme} point.  Straszewicz's theorem (e.g. Rockafellar, 1997, p. 167) states that the \emph{exposed} points of a closed convex set form a dense subset of
\emph{extreme} points of this set.
Thus, $\conv f(\cM)$ has an exposed point. Since an exposed
 point $y$ cannot be written as a \emph{strict convex combination} of other points of
 the set, $y$ must lie in $f(\cM)$. Let $x$ be a point from  $f^{-1}(y)$. Since $f$ is locally-convex at $x$ and there exists
 a hyperplane $H$
through $y$ that has empty intersection with $f(\cM) \setminus y$, we conclude that
the map $f$ is strictly locally-convex at $x$.
\end{proof}

Hereafter $i:\cM \rightarrow \X^n$  stands for a locally-convex complete immersion of a
 connected $(n-1)$-manifold $\cM$ into $\X^n$, and $M$ denotes $i(\cM)$. If $\cU \subset \cM$ we often use  $i:\cU \rightarrow \X^n$ for the restriction of $i|_{\cU}$ of $i$ to $\cU$. While discussing immersions, it is important to remember that they need not be injective; for example, we do not really consider, say, a line $L$ on the surface $M$ (as a set), but rather the map $i:\cL \rightarrow L$, where $\cL$ is 1-submanifold of $\cM$ and $L=i(\cL)$. The same philosophy applies in the case of any geometric subobjects of $i:\cM \rightarrow \X^n$. In the case of points we use the shorthand $i(x)$ instead of more proper $i:x \rightarrow i(x)$. Furthermore, when for some $\p \in M$ it is absolutely clear from the context as to which point $x$ of  $i^{-1}\p$ is considered, we may refer to $i:x \rightarrow i(x)=\p$ simply as ``point $\p$".

By a \emph{subspace} of $\X^n$ we mean an \emph{affine} subspace (i.e. a subspace
defined by a system $A\x=\bb$) in the case of $\R^n$, and the intersection of $\SSS^n
\subset \R^{n+1}$ with a \emph{linear} subspace of $\R^{n+1}$ in the case of $\SSS^n$.
A \emph{line} in $\X^n$ is a maximal geodesic curve; a \emph{plane} in $\X^n$ is a subspace of dimension 2; a \emph{hyperplane} is a subspace of $\X^n$ of codimension one. We will often use \emph{$k$-subspace} (or $k$-plane) instead of $k$-dimensional subspace -- the same convention applies to $k$-submanifolds, etc.
Two subspaces of   $\X^n$ are called complimentary if the sum of their dimensions is $n$ and the dimension of their intersection
is $0$ ($\dim \emptyset =-1$). For a subspace $S \subset \R^n$ we use  $\overrightarrow{S}$ to denote the linear space $S-S$.

For any $S \subset \X^n$ we denote by $\dim S$  the dimension of a minimal subspace containing $S$, which is denoted by $\aff S$ when it is  unique.
 The dimension operator $\dim$ has the meaning of topological dimension when applied to subsets of an abstract manifold, like $\cM$, and that of affine dimension when applied to subsets of $\X^n$ (i.e. for $S \subset \X^n$ we have $\dim S \triangleq \dim \aff S$).

\begin{theorem} (Jonker-Norman)  \label{theorem:J-N} Let $i:\cM \rightarrow \R^n$ ($n \ge 3$) be a complete locally-convex immersion of a connected $(n-1)$-manifold.
Then for any $x \in \cM$ there is a \emph{unique} submanifold $\cD$ through $x$ such that
\par \noindent 1) $i(\cD)=\aff i(\cD)$,
\par \noindent 2) $i|_{\cD}$ is a homeomorphism,
\par \noindent 3) $\cD$ is maximal with respect to 1) and 2).
\par Furthermore, for any hyperplane $H$ through $x$, which is complementary to $\aff \cD$ the set $\cG \triangleq i^{-1}(\cM \cap H)$ is a submanifold of $\cM$ such that
\par \noindent a) $\cM=\cD \times_{\mathrm{Top}} \cG$,
\par \noindent b) $i|_{\cG}$ is a locally-convex immersion into $\aff (\cM \cap H)$ with at least one point of strict convexity,
\par \noindent c) if $\cD'$ and $\cG'$ are to $x' \in \cM$ what $\cD$ and $\cG$ are  to $x$, then $\cD' \cong_{\mathrm Top} \cD$, $\cG' \cong_{\mathrm Top} \cG$, and $i(\cG)$ is equivalent to $i(\cG')$ under the action of affine automorphisms of $i(\cM)$ that map $i(\cD)$ to itself.
\par Finally, if $i$ \emph{is not a convex embedding,} then $\dim \cG =1$.
\end{theorem}

The theorem of Jonker and Norman generalizes Van Heijenoort's theorem by
characterizing the case of non-convex locally-convex (complete and connected) immersions.
Any such immersion is an immersion onto the product of a locally-convex, but not convex,
plane curve and a complimentary affine subspace. Whenever we have a map $i$ that satisfies
Jonker-Norman's theorem we will talk about the \emph{locally-convex direct
decomposition} of $i$; we may also say that the immersion $i$ \emph{splits} into the
\emph{locally-convex direct product} of $i:\cD\rightarrow D$ and $i:\cG\rightarrow G$. When
$i(\cG)$ is chosen to be perpendicular to $i(\cD)$, we call the decomposition
\emph{orthogonal.}

\begin{figure}[h]
\begin{center}
\resizebox{!}{100pt}{\includegraphics[clip=true,keepaspectratio=false]{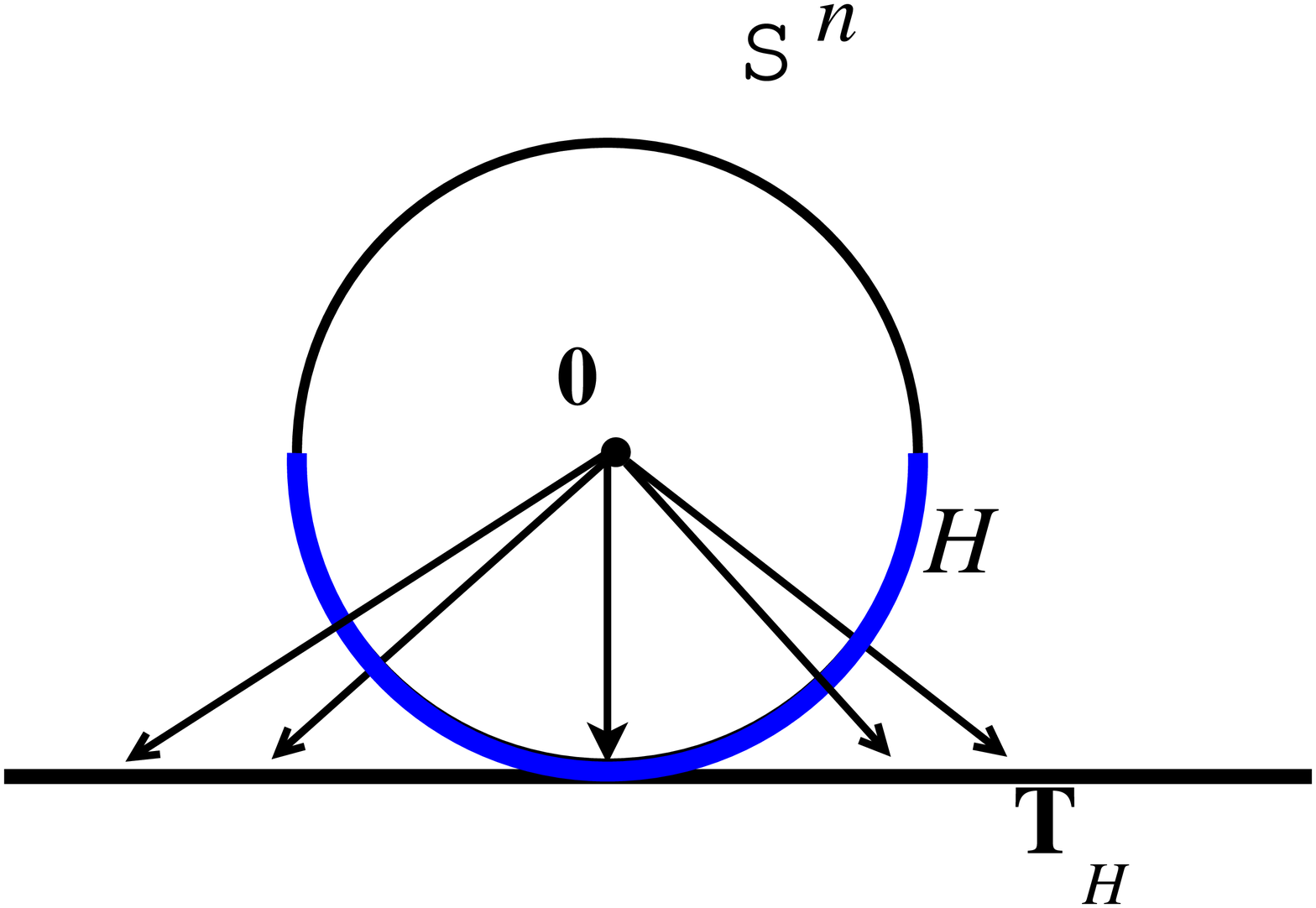}}
\caption{Central spherical projection map from hemisphere $H$ to $\bT_H$.} \label{fig:cp}
\end{center}
\end{figure}

By analogy with the traditional terminology for ruled surfaces and cylinders in 3D we refer
to $i:\cD \rightarrow D$, where $D= i(\cD)$, as a \emph{directrix} and to
$\overrightarrow{D}=D-D$ as the \emph{linear directrix} of $i:\cM \rightarrow \R^n$.
Similarly, we refer to  $i:\cG \rightarrow G$, where $G=i(\cG)$ , as a \emph{generatrix} of
$i:\cM \rightarrow \R^n$. Since the convex geometry of an (open) hemisphere of $\SSS^n$ is the same as that of $\R^n$, we will also use these terms on hemisphere (see Fig. \ref{fig:cp}).

\begin{figure}[h]
\begin{center}
\resizebox{!}{160pt}{\includegraphics[clip=false,keepaspectratio=false]{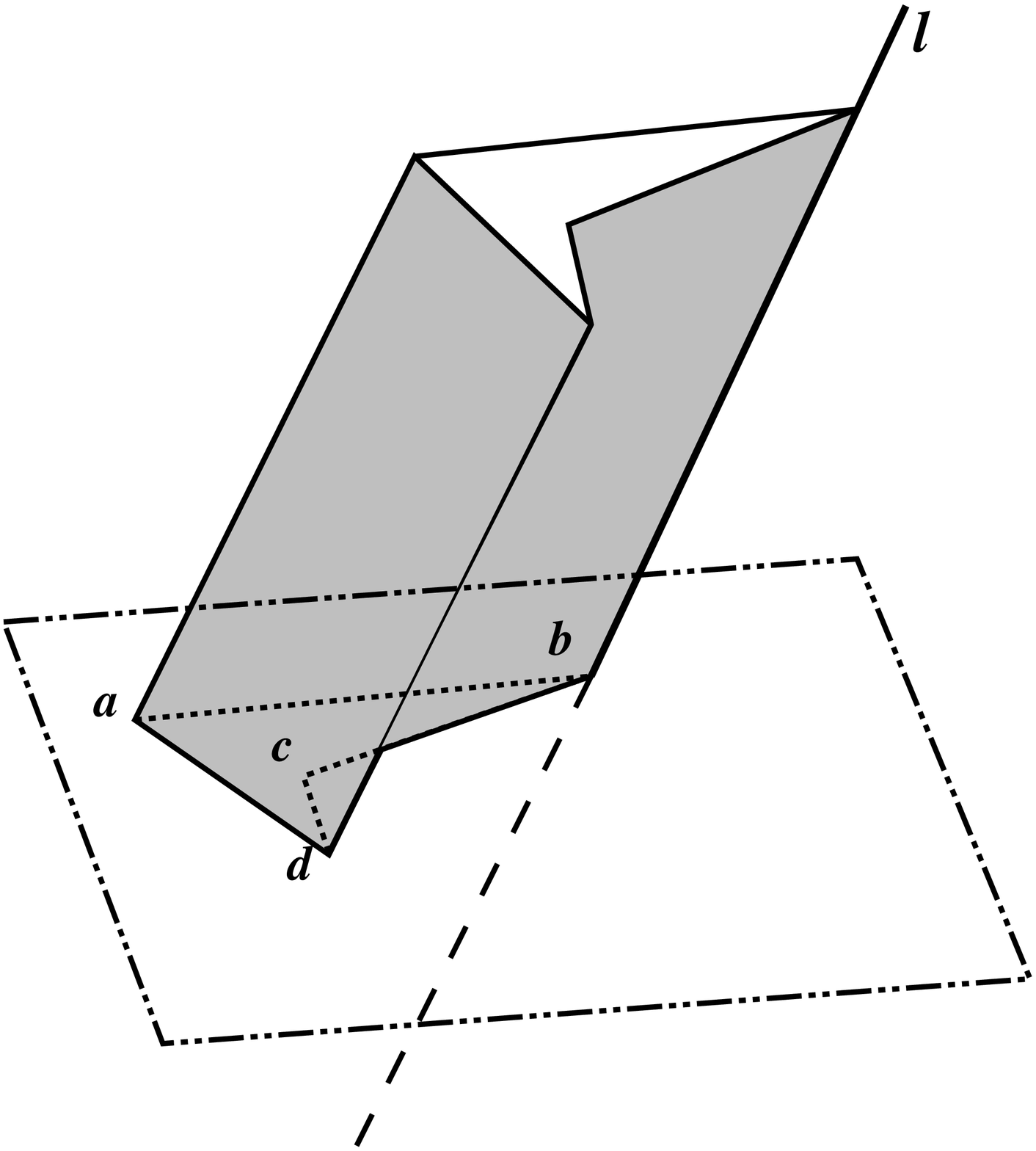}}
\caption{The product of a non-convex 4-gone and a line is locally convex, but not
globally convex} \label{fig:4-gone}
\end{center}
\end{figure}

By the Jonker-Norman theorem $\cD$ can be chosen so that it passes through a point of
strict convexity of $i\:|_{\cG}$, in which case we call it an \emph{exposed directrix.}
Note that such a directrix  is  \emph{never interior} to any flat of \emph{higher
dimension} contained in $i:\cM \rightarrow M=i(\cM)$. When needed we refer to $D$ as the geometric directrix and to $\cD$ as an abstract directrix, etc.  In Figure \ref{fig:4-gone}, which shows a piece of infinite cylindrical surface (it is embedded, so $\cM=M$ ), the line $\bl$ is a directrix and the 4-gone $(abcd)$ is a generatrix. There are exactly 4 exposed directrices and $a,b,c,d$ are the only points of strict convexity for the  section of the surface that is shown on the picture.

We call a connected submanifold $\cS$ of a topological space $\cT$ \emph{flat} with respect
to a realization map $r:\cT \rightarrow \X^n$ if $r: \cS \rightarrow r(\cS)$ is a
homeomorphism onto a subspace of $\X^n$  of the same dimension as $\cS$;   in this case we
call $r(\cS)$ a \emph{flat} contained in $(\cT,r)$. We will need the
following corollary of Theorem \ref{theorem:J-N}.

\begin{corollary}\label{cor:exposed directrix} In the context of Jonker-Norman theorem any
flat containing an exposed point of $i\:|_{\cG}$ is contained in  the exposed directrix
through this point.
\end{corollary}

The spherical convexity criterion, Theorem \ref{theorem:Spherical_Lemma}, is a direct
consequence of Theorem \ref{theorem:strict or conical} and Theorem \ref{theorem:union}; the
former deals with the case where a point of strict convexity is absent and the latter deals
with the case when it exists. The idea of proof of Theorem \ref{theorem:strict or conical}
is to apply Jonker-Norman's theorem locally, i.e. for a finite number of \emph{open}
hemispheres covering $\SSS^n$. The hypersurface, considered over each such hemisphere has a
number (possibly 0) of connected components, each of which having a unique orthogonal
Jonker-Norman decomposition (since the affine geometry of a hemisphere is essentially
equivalent to the geometry of $\R^n$). Among all such connected pieces of $(\cM,i)$ lying
in different hemispheres we pick one that has the exposed directrix of minimal dimension.
The Jonker-Norman decomposition is so ``rigid" that whenever an exposed directrix continues
from a hemisphere $H$ to a hemisphere $H'$ ($H \cap H' \neq \emptyset$), the Jonker-Norman
decompositions on $H \cap H'$, that are inherited from $H$ and $H'$,   must agree. As a result
we get an analog of Jonker-Norman's theorem for the sphere. Theorem \ref{theorem:union} is
proven by a combination of topological considerations and a metric (perturbation type)
argument reducing the problem to the Euclidean one.

\begin{theorem}\label{theorem:strict or conical}
Let $i: \cM \rightarrow \SSS^n$ ($n\ge3$) be a complete locally-convex immersion
of a connected $(n-1)$-manifold $\cM$ without points of strict convexity. Then
$i(\cM) = \mathbb{S}^n \cap \partial K$, where $K$ is a convex cone in
$\mathbb{R}^{n+1}$ containing the origin.
\end{theorem}

For a hemisphere $H \subset \SSS^n \subset \R^{n+1}$  we denote by $c_H$ the central
spherical projection map  (see Figure \ref{fig:cp})
 from $H$  onto the tangent $n$-plane
 $\mathbf{T}_{H}$ to $\SSS^n$ at the center of $H$ (when we find it convenient to index the tangent 
 space and the projection map by the center $\p$ of $H$ we write $\bT_{\p}$ and $c_{\p}$
 instead of $\bT_H$ and $c_H$).

 $\bT_H$ is an \emph{affine} real $n$-space; when we need to treat it as a \emph{linear} space (i.e.
$\overrightarrow{\bT_H}=\bT_H-\bT_H$), we identify the origin of  $\overrightarrow{\bT_H}$
with the center of $H$. First, let us make the following trivial observation.

\begin{lemma}\label{lem:extr-exp-sphere}
Let $\p$ be an \emph{extreme} point of a convex set $B \subset \SSS^n$. Then every
neighborhood of $\p$ has an \emph{exposed} point of $B$. In particular, if $\p$ is not
exposed itself, there are infinitely many exposed points on $\partial B$ arbitrarily close
to $\p$.
\end{lemma}
\begin{proof} This lemma is essentially a restatement for spherical spaces of a well-known theorem
 of Straszewicz on convex sets in $\R^n$ (Rockafellar, p.167). Since our lemma is
entirely of local nature, the proof in Rockafellar (1997) applies without changes.
\par Alternatively, consider the tangent space $\bT_{\p}$ to $\SSS^n \subset \R^{n+1}$ at $\p$.
The central projection maps the hemisphere $H$ centered at $\p$ onto $\bT_{\p}$ and $H \cap
B$ onto a convex set $B_{\p}$ in $\bT_{\p}$. The Euclidean theorem can now be applied to
$\p$ as an extreme point of the convex set $B_{\p}$ in $\bT_{\p}\cong \R^n$. The property
of a point to be extreme or exposed is preserved under the central projection and its
inverse.
\end{proof}

\medskip

\begin{proof}[Proof of Theorem \ref{theorem:strict or conical}]
If there is $x \in \cM$ such that  $i(x)$ is extreme for some  convex witness at $x$, then, by Lemma
\ref{lem:extr-exp-sphere}, there is an exposed point at every neighborhood of $x$. Since an exposed
point of a convex body is the same as a point of strict convexity of its surface,  $(\cM,i)$ is strictly
convex in at least one point, which is impossible. Thus, for each $x \in \cM$  the image $i(x)$ is an interior point of a segment on $M \triangleq i(\cM)$.

For an open hemisphere $H$ let $\cSH$ denote the set of all maximal connected submanifolds
of $i^{-1} (H)$  whose $i$-images lie in $H$. When $\cSH\neq \emptyset$ each $\cS\in \cSH$
is an $(n-1)$-submanifold of $\cM$ and $c_H \circ i:\cS \rightarrow \bT_H$ is a complete
locally-convex immersion (e.g. by Lemmas \ref{lem:complete proper} and \ref{lem:proper
complete}).

By Jonker-Norman's theorem $c_H \circ i :\cS \rightarrow \bT_H$ has a
 \emph{locally-convex direct decomposition} $\cS=\cG \times \cL$, where
  $c_H \circ i\:|_{\cG}$ is a locally-convex immersion  of a compact connected
  $g$-submanifold $\cG$ and $c_H \circ i:\cL\rightarrow L$
  is a homeomorphism from a $d$-submanifold $\cL \subset \cS$ onto a $d$-subspace of $\bT_H$ ($n-1=d+g$).
  Denote by $\vL(\cS)$ the linear space $L-L$.

Let us pick $\cG$ so that  $G \perp L$ where $G \triangleq c_H \circ i(\cG)$. Then $c_H
\circ i\:|_{\cS}$ is the
 \emph{orthogonal direct product} of the generatrix $c_H \circ i\:|_{\cG}$  and the directrix
 $c_H \circ i\:|_{\cD}$.  On the hemisphere $H$ this decomposition corresponds to the orthogonal locally-convex split of $i\:|_{\cS}$
into a hemispherical generatrix $i: \cG \rightarrow c_H^{-1}G$  and a hemispherical
directrix $i:\cD \rightarrow c_H^{-1}L$.

Let $\mfH$ be a finite covering of $\SSS^n$ by \emph{open hemispheres}. Since $i$ does not
have points of strict convexity,  $i(\cM)$ is not completely covered by any  single
hemisphere. We will use  $\cS_H$ for an element of $\cSH$ -- subindex $H$ only indicates
that $\cS_H$ was chosen from $\cSH$.
 Likewise, once $\cS_H$ is fixed, we may use $L_{\cS_H}$ and $\cG_{\cS_H}$, etc.
 to indicate that $L_{\cS_H}$ and $\cG_{\cS_H}$ are  obtained from the direct decomposition of
 $c_H\circ i\:|_{\cS_H}$.

Suppose  $U=i({\cU})$ is a convex hypersurface for some connected submanifold $\cU \subset
\cS$. Let $H,H' \in \mfH$, and $H \cap H' \cap U \neq \emptyset$. Then there is a unique
$\cS_{H'} \in \cSH'$ such that $\cS_{H'}$ contains $i^{-1}H' \cap \cU$. We will refer to
this fact by saying that \emph{whenever a convex subsurface of $i:\cS_{H}\rightarrow H$
protrudes into $H'$, the surface $i:\cS_{H}\rightarrow H$ extends uniquely into $H' \cup H$
along $U$} (or, in other words, the map $i\:|_{\cS_{H}}$ extends uniquely over $i^{-1}H'$
along $\cU \cap i^{-1}H'$) . In this context $\cS_{H'} \cup \cS_H$ is called the
\emph{extension} of $\cS_{H}$ and $\cS_{H'}$ is called an \emph{adjoint} to $\cS_{H}$.

Among the elements of $\mfH$ that overlap with $i(\cM)$, let $H_0$ be one where we can pick
$\cS_{H_0} \in \cSH_0$ so that $d \triangleq\dim \vL(\cS_{H_0})\le \dim \vL(\cS_H)$ for all
$H$ that overlap with $i(\cM)$ and each possible choice of $\cS_H \in \cSH$; let $i:\cD_0
\rightarrow D_0=i(\cD_0)$ be an exposed hemispherical directrix for $\cS_{H_0}$.

\medskip
If $d=0$, then $\cS=\cG$, where $c_{H}\circ i:\cG \rightarrow \bT_H$ has a point of strict
convexity, which contradicts to our assumption about $i$.

\medskip
If $d=n-1$, then $\aff D_0$ is an $(n-1)$-hemisphere of $H_0$ and $\cS_{H_0}=\cD_0$. We
know that whenever a convex subsurface of $i:\cS_{H_0}\rightarrow H_0$ protrudes into $H$,
the surface $i:\cS_{H_0}\rightarrow H_0$ extends uniquely into $H$ along this subsurface,
which implies that $i\:|_{\cD_0}$ can be extended to all hemispheres overlapping with $\aff
D_0$. Since $\cM$ is a connected $(n-1)$-manifold, $i:\cM \rightarrow \SSS^n$ is an
immersion onto $\aff \cD_0$, which is, by the covering mapping (see Seifert \& Threlfall)
theorem, a homeomorphism if $n>2$.

\medskip
Let $1 \le d \le n-2$.  Let $H\cap D_0 \neq \emptyset$ for some $H \in \mfH$.  We know
$i\:|_{\cS_{H_0}}$ extends in a unique way along $D_0 \cap H$ into $H_0 \cup H$: denote the
adjoint element of $\cSH$ by $\cS_{H}$. Obviously, $\aff D_0 \cap H$ is an extreme
(geometric) hemispherical  directrix for $i\:|_{\cS_{H}}$. As $D_0$ is completely covered
by elements of $\mfH$, the submanifold $\cD_0$ extends to a connected component of
$i^{-1}(\aff D_0)$   inside of $\cM$. Set $D\triangleq \aff D_0$ and let  $\cD$ be a
maximal connected $d$-submanifold  of $\cM$ such that $D=i(\cD)$.  Since $i$ is a complete
immersion, it is proper (preimages of  compact sets are compact) and $\cD$ is compact.
Thus, the preimage of any $\p \in D$ under  $i\:|_{\cD}$ is a finite set of the size that
does not depend on $\p$.

Without loss of generality we  assume that $D$ is completely covered by hemispheres
\newline $H_0,\dots,H_N \in \mfH$,~all centered at points of $D$; denote this subset of
$\mfH$ by $\mfH_D$. Let $\cS$ be a connected component of $i^{-1}(\cup_{j=0}^{N}H_j)$ that
contains $\cD$, i.e. $\cS$ is the unique maximal extension of $\cS_{H_0}$ into
$\cup_{j=0}^{N}H_j$ along $D$. For $H_k,H_l \in \mfH_D$, where $H_k \cap H_l \neq
\emptyset$, on each connected component of $i^{-1}(H_k \cap H_l) \cap \cS$ the
locally-convex orthogonal decompositions of $i\:|_{\cS}$, which are induced by   the
restrictions of $i$ to $\cS \cap i^{-1}H_k$ and $\cS \cap i^{-1}H_l$ respectively, agree;
this follows directly from Jonker-Norman's theorem. Furthermore, since both $H_k$ and $H_l$
are centered at $D$, the generatrices in these two locally-convex \emph{orthogonal}
decompositions are all isometric to each other -- the rotational subgroup $Iso^+(D)$ of
$Iso(D)$ is transitive on them.

 Thus, we have a \emph{locally-convex orthogonal  fibration} of the immersion $i\:|_{\cS}$:
namely, we have a continuous map $\pi:\cS \rightarrow \cD$, which sends each (topological)
generatrix into its base point on $\cD$ and for each $x \in \cD$ there is a neighborhood
$\cU_x \subset \cD$ such that $\pi^{-1}(\cU_x)$ is the direct orthogonal product of $\cU_x$
and a fiber $\cG_x$ over $x$, such that $i\:|_{\cG_x}$ is a locally-convex immersion into
$D^{\perp}_x$, where the latter is the orthogonal complement of $D$ through $x$. Inside of
each $H \in \mfH_D$ the fibers (i.e. generatrices) are isometric; moreover, as we just
noticed above, the fibers from different $H$'s are also isometric. Thus, the constructed
locally-convex fibration of the immersion $i\:|_{\cS}$ is, in fact, a direct product
decomposition, i.e. $\cS=\cD \times \cG$, where $\dim \cG=n-1-d$ and $i\:|_{\cG}$ is a
locally convex immersion into a $(n-d)$-hemisphere perpendicular to $D$.

Set $D^*\triangleq \SSS^n \cap (\cone D)^{\perp}$, where $\cone D$ is the cone with apex at
$\mathbf{0}$ over $D$. $D^{*}$ consist of all points of $\SSS^n$ that are not covered by
the elements of $\mfH_D$. We claim that all generatrices from the orthogonal decomposition
of $i\:|_{\cS}$ ``reach to $D^*$", i.e. for each generatrix $\cG \subset \cS$ and any
neighborhood of $D^*$ there is $p \in \cG$ such that $i(p)$ lies in this neighborhood.   By
contradiction: let $p \in \cG \subset \cS$ be such that the distance $\rho>0$ between
$i(p)$ and $D^{*}$ is equal to the distance between $i(\overline{\cG})$ and  $D^{*}$. Since
all generatrices are isometric with respect to $Iso^+(D)$, $\cS$ contains a submanifold
mapped onto the orbit of $p$ under the induced action of $Iso^+(D)$ on $\cS$, but does not
contain any points mapped by $i$ to spherical points at the distance smaller  than $\rho$
from $D^{*}$. Then $i$ is not locally-convex at all points of this submanifold. Thus,  all
generatrices of the orthogonal decomposition of $i\:|_{\cS}$ ``reach to $D^*$".

Since $\cM$ is compact, for any  $p \in \overline{\cS} \setminus \cS$ we have $i(p) \in
D^{*}\cong \SSS^{n-d-1}$. Then $i(p)$ belongs to the closure of each  generatrix from the
orthogonal  fibration of $i\:|_{\cS}$ with base $i\:|_{\cD}$. Let $H(\p)$ be the hemisphere
centered at $\p=i(p)$. Under $c_{H(\p)}:H_p \rightarrow \bT_{H(\p)}$ the points of  $D$
correspond to rays emanating from the origin of $\bT_{H(\p)}$, or, in other words, in ``the
world of" $\bT_p$ the spherical subspace $D$ corresponds to a ``$d$-sphere at infinity",
which we denote by $D(\bT_{H(\p)})$. Thus, the isometry group of the surface
$(\overline{\cS}, c_{H(\p)}\circ i)$ includes all linear isometries (in particular,
rotations about $\p$) that preserve the sphere $D(\bT_{H(\p)})$ at infinity. We know that
$\p=i(p)=c_{H(\p)}\circ i(p)$ must belong to the interior of a  segment $I=c_{H(\p)}\circ
i(\cI)$ on this surface. Any isometry of $(\overline{\cS}, c_{H(\p)}\circ i)$ will map $I$
to another line segment. Because of local convexity at $p$, the isometries preserving the
sphere $D(\bT_{H(\p)})$ at infinity belong to the isometry group of a supporting hyperplane
at $\p=c_{H(\p)}\circ i(p)$. Since $I$ must be in this hyperplane, $i(p)$ is interior to a
$(d+1)$-flat of $\overline{\cS}$. We will have to deal separately with the cases $d=n-2$
and $d \le n-3$.

\textbf{Case: $d=n-2$.} $\overline{\cS} \setminus \cS \cong \SSS^0$. Then $i(p)$ is
interior to an $(n-1)$-flat $i:\cF \rightarrow F$ of $(\overline{\cS},i)$. We will show
that $i\:\:|_{H(\p)}$ is a homeomorphism onto an open $(n-1)$-hemisphere of $\SSS^n$.
Consider a directed geodesic in $\aff F$ with source at $i(p)$  that does not extend to $D$
inside of $(\overline{\cS},i)$. Let $\bb=i(b)$ be a point where this geodesic first
diverges with the surface $(\overline{\cS},i)$. Point $b$ belongs to a unique fiber from
the orthogonal fibration of $(\overline{\cS},i)$ with base $(\cD,i)$. The isometry group of
$D$ is transitive on the fibers. Thus, all directed geodesics through $i(p)$ that lie in
$(\overline{\cS},i)$ diverge with the $(n-1)$-flat $F$ at the same distance from $i(p)$;
hence, $\cF$ is an $(n-1)$-ball (in the $i$-distance) centered at $p$. But then all points
of $\partial \cF$ are extreme points for $(\overline{\cS},i)$, which contradicts our
assumptions. Thus $F$ is an $(n-1)$-hemisphere centered at $i(p)$ and bounded by $D$. The
same argument is applied to the other point of $\overline{\cS} \setminus \cS \cong \SSS^0$.
Thus, $i$ is a homeomorphism onto the surface made of two $(n-1)$-hemispheres glued
together at their common $(n-2)$-dimensional boundary $D$.

\textbf{Case: $1 \le d < n-2$.} The central projection of a generatrix $i\:|_{\cG}$ (where
$G=i(\cG)$ and  $G\perp D$) with a base point $i(x)=\x \in D$ onto its tangent subspace
$\bT_G \subset \bT_{\x}$ at $\x \in D$ is a locally-convex unbounded complete surface
$c_{\x} \circ i\:|_{\cG}$. Since the topological dimension of the generatrix is larger than
one, by Jonker-Norman's theorem  it is an embedded convex surface in $\bT_G$ and $\x$ is a
(geometric) point of strict convexity for this surface.  Thus, $i\:|_{\cG}$ is a convex
surface on $H_{\x}$. We need to understand the geometry of $i\:|_{\cG}$ at infinity, i.e.
at $\partial H_{\x} \cap D^*$. Because of strict convexity at $\x=i(x)$,  $\partial H_{\x}
\cap i(\overline{\cG})$ is the boundary of a strictly convex compact set in  $\partial
H_{\x}$. Suppose $z \in \cG$ is a point of strict convexity of $i\:|_{\cG}$ and $z \neq x$.
Then there is an exposed hemispherical directrix $i:\cD_z \rightarrow H_p$ through $z$,
distinct from $\cD \cap i^{-1}H_{\x}$. Since $i(\cD)$ and $i(\cD_z)$ are parallel in
$H_{\x}$, they ``intersect at infinity" (i.e. on $\partial H_{\x}$) over a common
$(d-1)$-sphere. Thus, we have \emph{two distinct exposed directrices} through the same
point of $\cM$. This is impossible by Corollary \ref{cor:exposed directrix}. Thus,
$i\:|_{\cG}$ has a unique point of strict convexity.

$c_{\x} \circ i:{\cG}\rightarrow \bT_G$ is an embedded unbounded complete convex
hypersurface with a unique point of strict convexity. Let us apply a projective
transformation that sends $\bT_G$ to another subspace $\bP$ of the same dimension in
$\R^{n+1}$ in such a way that the point $\x$ of $G \subset  \bT_G$ is mapped to a point at
infinity of $\bP$. This will give us an embedded unbounded   complete convex hypersurface
in $\bP$ without points of strict convexity. By   Jonker-Norman this hypersurface in $\bP$
is the product of a line $L$ in $\bP$ and a   compact convex hypersurface in a subspace of
$\bP$, which is complimentary to $L$. Thus,   $c_{\x} \circ i ({\cG}) $ is the boundary of
a cone with apex at $\x$ over a convex   compact set ``on the sphere at infinity of
$\bT_G$". Hence, $c_{\x}\circ i\:|_{\overline{\cG}}$ is an immersion onto a cone over a
convex compact surface of topological dimension $(n-1)-d-1=n-d-2$ on $D^*$.   Since all
generatrices are isometric to $i\:|_{\overline{\cG}}$ with respect to the action of
$Iso^+(D)$,  we conclude that $\cM$ contains a closed  $(n-1)$-submanifold $\overline{\cS}$
(without boundary). Since $\cM$ is connected, $\overline{\cS}=\cM$.

\textbf{Remark on injectivity:} the proof does not imply that $i$ is an embedding. When there are no points of strict convexity, non-injectivity is possible if and only if $d=\dim \cD=1$. When $d>1$ the classical covering mapping theorem (see e.g. Seifert-Threlfall)  implies that the map is one-to-one.
\end{proof}
\medskip

When a minimal geodesic between $\p$ and $\q$ is unique, we denote it by $[\p,\q]$; we will also use $[p,q]$, where $i(p)=\p$, $i(q)=\q$, to refer to a curve in $\cM$ that is mapped homeomorphically onto  $[\p,\q]$.
\begin{theorem} \label{theorem:union}Let $\cM$ be connected and
let $i:\cM \rightarrow \X^n$ be complete, locally-convex and, also, strictly
locally-convex at $o \in \cM$.
Then $i:\cM \rightarrow \X^n$ is a convex embedding.
\end{theorem}

\begin{proof} If $\bH_o \subset \SSS^n$ is a supporting hyperplane at $i(o)$, then
let us denote the \emph{open} hemisphere defined by $\bH_o$ that contains the image of a small neighborhood of $o$ by $\bH^+_o$; the other open hemisphere is then  denoted by $\bH_o^-$. If $\cN$ is a neighborhood of $x$ we denote by $\dot{\cN}$  its punctured
version, i.e. $\cN \setminus x$.

Let $\cS$ be a maximal connected $(n-1)$-submanifold of $\cM$ such that $o \in \cS$ and $i(\dot{\cS}) \subset \bH_o^+$. Suppose that there is
no $x \in \overline{\cS} \setminus o$ with $i(x) \in \bH_o$. Then the distance between
$\overline{\cS} \setminus \cN_o$ (where $\cN_o$ is a small neighborhood of $o$)  and
$\bH_o$ is strictly positive. This means we can perturb $\bH_o$ so that $i(\overline{\cC_{o}})$ is in $\widetilde{\bH}^+$, where $\widetilde{\bH}$ is a perturbed version of $\bH_o$.  Let $c$ be the
central projection on $\bT_{\widetilde{\bH}^+}$ from $\widetilde{\bH}^+$. Clearly, $c \circ i\:|_{\cS}$ satisfies the
conditions of Van Heijenoort's theorem; hence, $c \circ i\:|_{\cS}$ is a convex embedding.
Since $\cM$ is connected,  $\cS=\cM$.

Let now $p \in \overline{\cS}\setminus o$ be such that  $i(p) \in \bH_o$. If $i(p)
\neq i(o)^{\mathrm{op}}$ ($i(o)^{\mathrm{op}}$ stands for  the opposite of $i(o)$), then the minimal geodesic joining $i(o)$ and $i(p)$ is unique and lies in $\bH_o$.
Let $\{i:[o,x_m] \rightarrow [i(o),i(x_m)]\}_{m \in \N}$, with  $[o,x_m] \subset \cM$, be a sequence of minimal geodesics that converges to $i:[o,p] \rightarrow [i(o),i(p)]$. The geodesics in this sequence lie arbitrarily close to $\bH_o$. Since $(\cM,i)$ is strictly convex at $o$, we find that  $p=o$, which contradicts to the choice of $p$. Thus,
the points of $\overline{ \cS}\setminus o$ that are mapped to $\bH_o$ are mapped to   $i(o)^{\mathrm{op}}$. Since $i$ is a proper immersion, the preimage of $i(o)^{\mathrm{op}}$ in $\cM$ is finite. Hence, the preimage of $i(o)^{\mathrm{op}}$ in $\partial \cS=\overline{ \cS}\setminus \cS$ is finite.
 Clearly, $c \circ i\:|_{\dot{\cS}}$
satisfies the conditions of Jonker-Norman's theorem. Since $i$ is strictly convex at
$o$, $c \circ i\:|_{\cS}$ must be a convex unbounded \emph{embedding} onto a cylinder in
$\bT_{\bH_o^+}$. The directrix must be 1-dimensional, for the cylinder has only two points at infinity, $i(o)$ and $i(o)^{\mathrm{op}}$ (see Fig. \ref{fig:dolya}). Thus, $\cS$ contains a punctured neighborhood of $p$ homeomorphic to an $(n-1)$-ball. If we add $p$ to $\cS$ we get a compact connected $(n-1)$-submanifold of $\cM$ (without boundary). Since $\cM$ is connected, $\overline{\cS}=\cM$.
\begin{figure}[h]
\begin{center}
\resizebox{!}{300pt}{\includegraphics[clip=true,keepaspectratio=false]{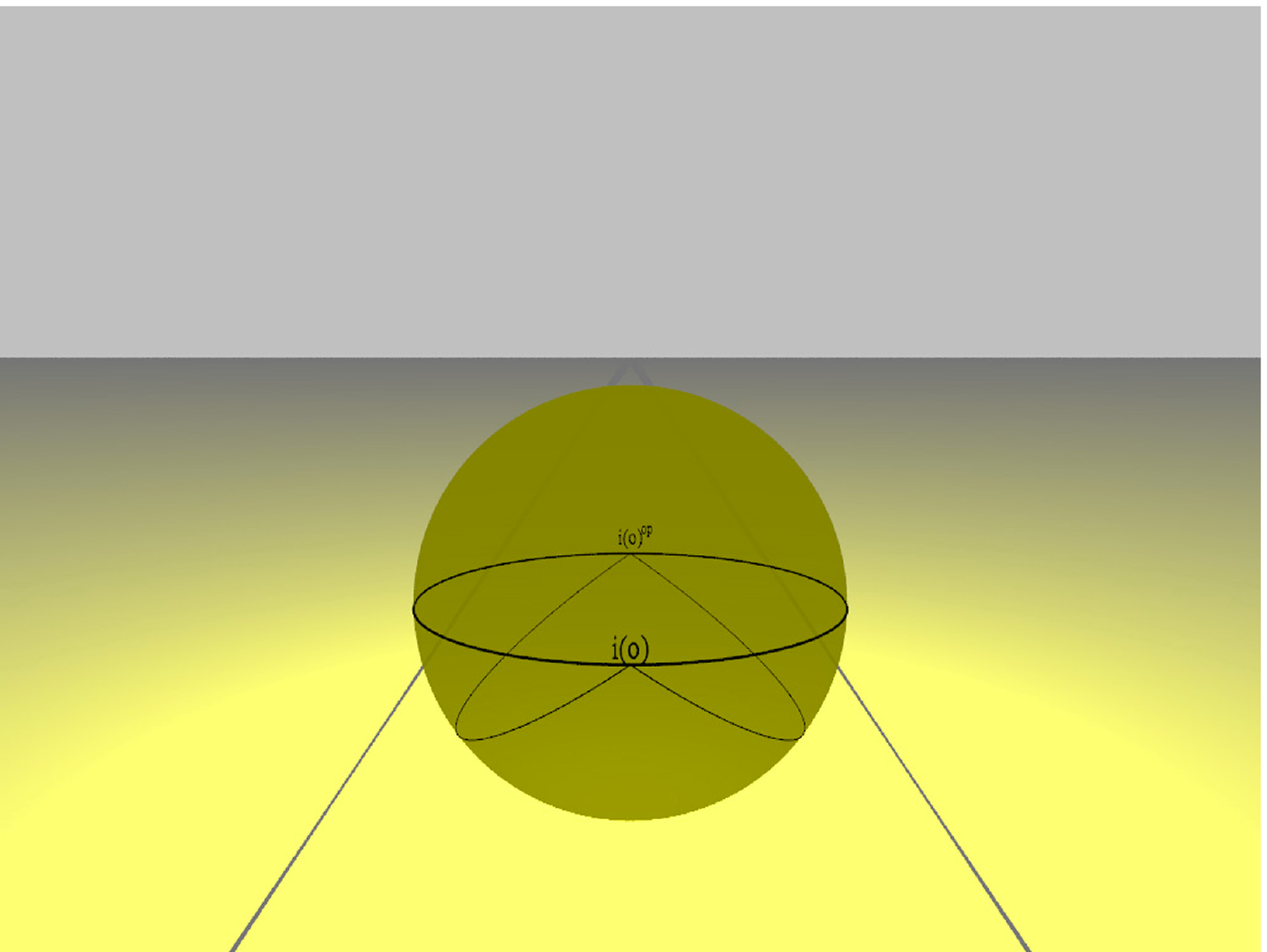}}
\caption{The surface $i\:|_{\overline{\cS}}$ has only two points on $\bH_o$: $i(o)$ and $i(o)^{\mathrm{op}}$ } \label{fig:dolya}
\end{center}
\end{figure}
Thus, $i:\cM \rightarrow \X^n$ is a convex embedding.
\end{proof}

\medskip
Konstantin Rybnikov
\par \url{Konstantin_Rybnikov@uml.edu}
\par University of Massachusetts at Lowell, Mathematical Sciences, Lowell, MA 01854 USA
\end{document}